\pdfoutput=1
\nonstopmode
\documentclass{amsart}

\usepackage{amsmath, amssymb, amsthm, amscd}
\usepackage{calc}
\usepackage{psfrag,graphicx}
\usepackage[notightpage]{pst-pdf}
\usepackage{verbatim}
\usepackage[all]{xy}
\newdir{ >}{{}*!/-10pt/@{>}}
\newdir^{ (}{{}*!/-5pt/@^{(}}
\newdir_{ (}{{}*!/-5pt/@_{(}}
\def\cgaps#1{}
\def\Cgaps#1{}
\allowdisplaybreaks

\def\undersetbrace#1\to#2{\underbrace{#2}_{#1}}
\def\oversetbrace#1\to#2{\overbrace{#2}^{#1}}
\def\AMSunderset#1\to#2{\underset{#1}{#2}}
\def\AMSoverset#1\to#2{\overset{#1}{#2}}

\def\norm#1{\left\|{#1}\right\|}

\newtheorem*{prop*}{Proposition}

\newtheorem*{thm*}{Theorem}

\newtheorem*{lem*}{Lemma}

\newtheorem*{cor*}{Corollary}
\numberwithin{equation}{subsection}

\def\ign#1{}             
\def\o{\circ}
\def\X{\mathfrak X}

\def\ka{\kappa}

\def\Ga{\Gamma}

\def\Ph{\Phi}

\def\i{^{-1}}
\def\x{\times}
\def\p{\partial}
\let\on=\operatorname
\def\L{\mathcal L}
\def\ol{\overline}
\def\AMSonly#1{}

\def\R{\mathbb{R}}
\def\ol{\overline}
\def\Tr{\on{Tr}}
\def\vol{\on{vol}}
\def\Imm{\on{Imm}}
\def\g{\bar{g}}

\def\Vol{\on{Vol}}

\begin{document}
\title{Curvature weighted metrics on shape space of hypersurfaces in $n$-space}
\author{Martin Bauer, Philipp Harms, Peter W. Michor}
\address{
Martin Bauer, Philipp Harms, Peter W. Michor:
Fakult\"at f\"ur Mathematik, Universit\"at Wien,
Nordbergstrasse 15, A-1090 Wien, Austria.
}
\email{Bauer.Martin@univie.ac.at}
\email{Philipp.Harms@univie.ac.at}
\email{Peter.Michor@univie.ac.at}
\thanks{All authors were supported by 
FWF Project 21030 and by the 
    NSF-Focused Research Group: 
    The geometry, mechanics, and statistics of the infinite dimensional
    shape manifolds.}
\date{\today}
\begin{abstract}
Let $M$ be a compact connected oriented $n-1$ dimensional manifold without boundary.
In this work, shape space is the orbifold of unparametrized immersions from $M$ to $\mathbb R^n$.
The results of \cite{Michor118}, where mean curvature weighted metrics were studied, 
suggest incorporating Gau{\ss} curvature weights in the definition of the metric. 
This leads us to study metrics on shape space that are induced by metrics on the space of immersions
 of the form 
$$ G_f(h,k) = \int_{M} \Ph . \g(h, k) \on{vol}(f^*\bar{g}).$$
Here $f \in \Imm(M,\R^n)$ is an immersion of $M$ into $\R^n$ and $h,k\in C^\infty(M,\mathbb R^n)$ 
are tangent vectors at $f$. 
$\bar g$ is the standard metric on $\mathbb R^n$, $f^*\bar g$ is the induced metric on $M$,
$\on{vol}(f^*\bar g)$ is the induced volume density and 
$\Ph$ is a suitable smooth function depending on the mean curvature and Gau{\ss} curvature.
For these metrics we compute the geodesic equations both on the space of immersions and on shape 
space and the conserved momenta arising from the obvious symmetries.
Numerical experiments illustrate the behavior of these metrics. 
\end{abstract}

\maketitle


\section{Introduction}\label{sec:introduction}

Riemannian metrics on shape space are of great interest for many procedures in data analysis
and computer vision. They lead to geodesics, to curvature and diffusion. 
Eventually one also needs statistics on shape space. 
The shape space used in this paper is the orbifold $B_i(M,\R^n)$ of unparametrized immersions of
a compact oriented $n-1$-dimensional manifold $M$ into $\R^n$, see \cite{Michor40}:
$$ B_i(M,\R^n) = \Imm(M,\R^n) / \on{Diff}(M). $$
The shape space $B_e$ of submanifolds in $\R^n$ of type $M$ used in \cite{Michor118,Michor119} is a special case of this. 

The quest for suitable Riemannian metrics was started by Michor and Mumford in \cite{Michor98,Michor102}, where
it was found that the simplest and most natural such metric induces vanishing geodesic distance on shape space $B_i$.
One attempt to strengthen the metric is by adding weights depending on the curvature. 
For curves, where there is only one notion of curvature, this was done in \cite{Michor98,Michor107}. 
In higher dimensions, there are several curvature invariants. We are especially interested in 
the case $n=3$, where the curvature invariants are the mean and Gau{\ss} curvature. 
Mean curvature weights were used in \cite{Michor118}. 
In this paper we take up these investigations
by adding a Gau{\ss} curvature term to the definition of the metric. 

A Riemannian metric on $\on{Imm}(M,\mathbb R^n)$ is a family of
positive definite inner products $G_f(h,k)$ where $f \in
\on{Imm}(M,\mathbb R^n)$ and $h,k \in C^\infty(M,\mathbb R^n) = T_f\on{Imm}(M,\mathbb R^n)$
represent vector fields on $\mathbb R^n$ along $f$. 
We require that our
metrics will be invariant under the action of $\on{Diff}(M)$, hence the
quotient map resulting from dividing by this action will be a Riemannian submersion. 
This means that the tangent map of the quotient map $\Imm(M,\R^n)\to B_i(M,\mathbb R^n)$ is 
a metric quotient mapping between all tangent spaces.
Thus we will get Riemannian metrics on $B_i$. 
Riemannian submersions have a very nice effect on geodesics: the geodesics on
the quotient space $B_i$ are exactly the images of the horizontal geodesics on the top
space $\on{Imm}$. 

The metrics we will look at are of the form:
$$ G_f(h,k) = \int_{M} \Ph\big(\on{Tr}(L),\on{det}(L)\big) \g( h, k ) \on{vol}(f^*\bar{g})$$
Here $\bar g$ is the standard metric on $\mathbb R^n$, $f^*\bar g$ is the induced metric on $M$,
$\on{vol}(f^*\bar g)$ is the induced volume density, $L$ is the Weingarten mapping of $f$ and 
$\Ph: \R^2 \to \R_{>0}$ is a suitable smooth function depending on the mean curvature
$\Tr(L)$ and Gau{\ss} curvature $\det(L)$.
Since $\Phi$ is just multiplication with a smooth positive function, 
the horizontal bundle consists exactly of those vector fields which are pointwise normal to 
$f(M)$, see \cite{Michor118}. 

Big parts of this paper can be found in Martin Bauer's Ph.D. thesis \cite{Bauer2010}.
\section{Differential geometry of surfaces and notation}\label{notation}

We use the notation of \cite{MichorH}. Some of the definitions can also be found in \cite{Kobayashi1996a}.
Similar expositions in the same notation are \cite{Michor118, Michor119}.

In all of this chapter, let $f:\R \to \on{Imm}(M,\R^n)$ be a smooth path of immersions. By convenient calculus
\cite{MichorG}, $f$ can equivalently be seen as a smooth mapping $f:\R \x M \to \R^n$ 
such that $f(t,\cdot)$ is an immersion for each $t$. 
We will deal with the bundles
\begin{equation*}\xymatrix{
T^r_s M \ar[d] & 
T^r_s M \otimes f(t,\cdot)^*T\R^n \ar[d] & 
\on{pr}_2^* T^r_s M \otimes f^*T\R^n \ar[d]\\
M & M & \R \x M
}\end{equation*}
Here $T^r_sM$ denotes the bundle of $\left(\begin{smallmatrix}r\\s\end{smallmatrix}\right)$-tensors on $M$, i.e.
$$T^r_sM=\bigotimes^r TM \otimes \bigotimes^s T^*M \quad \text{ and } \quad \on{pr}_2:\R \x M \to M.$$

Let $\bar g$ denote the Euclidean metric on $\mathbb R^n$. 
The \emph{metric} induced on $M$ by $f(t,\cdot)$ is the pullback metric $g=f(t,\cdot)^*\bar g$. 
On tensor spaces $T^r_s M$ we consider the product metric $g^r_s:=\bigotimes^r g \otimes \bigotimes^s g\i$.
A useful fact is (see \cite{Michor119})
\begin{equation}
g^0_2(B,C)= \on{Tr}(g\i B g\i C) \quad \text{for $B,C \in T^0_2M$ if $B$ or $C$ is symmetric.}
\end{equation}

Furthermore $g^r_s \otimes \bar g$ yields a metric on $T^r_s M \otimes f(t,\cdot)^*T\R^n$. 
Let $\nabla$ denote the \emph{Levi-Civita covariant derivative} acting
on sections of all of the above bundles.
The covariant derivative admits an adjoint, which is given by $\nabla^*=-\Tr^g \nabla$. For more details see \cite{Michor119}.

The \emph{normal bundle} of an immersion $f(t,\cdot)$ or a path of immersions $f$, respectively, 
is given by
\begin{equation*}\xymatrix{
\big(Tf(t,\cdot)\big)^\bot \ar[d] \ar @{} [rrd]|*{\text{or}} && Tf^\bot \ar[d] \\
M && \R \x M
}
\end{equation*}

Let $\nu$ designate the positively oriented unit normal field with respect to the orientations on $M$ and $\R^n$.
Any field $h$ along $f(t,\cdot)$ or $f$ can be decomposed uniquely as
$$h=Tf(t,\cdot).h^\top + h^\bot.\nu$$ 
into parts  which are tangential and normal to $f$. The two parts are defined by the relations
\begin{align*}
&h^\bot =\g( h,\nu )\in C^\infty(M) \\
&h^\top \in \X(M), \text{ such that } g(h^\top,X)=\g(h,Tf(t,\cdot).X) \text{ for  all } X\in \X(M).
\end{align*}

Let $s^f$ denote the \emph{second fundamental form} of $f(t,\cdot)$ on M and $L$ the 
\emph{Weingarten mapping}, which are related to the Levi-Civita covariant derivative 
$\nabla$ on $(\R^n,\g)$  by the following equation:
$$ \g(\nabla_X \nu, Tf(t,\cdot).Y ) = -s(X,Y)=- g( LX,Y ). $$

The eigenvalues of $L$ are called \emph{principal curvatures} and
the eigenvectors \emph{principal curvature directions}. 
$\on{Tr}(L)$ is the \emph{mean curvature}, 
and in dimension two $\det(L)$ is the \emph{Gau{\ss}-curvature} (we use this name generally).

\section{Variational formulas}\label{variation}

We will calculate the derivative of the Gau{\ss}-curvature. The differential calculus used is convenient calculus as in 
\cite{MichorG}.
Proofs of some of the formulas  in this chapter can be found in \cite{Besse2008,Michor102,Michor118,Michor119} and
in \cite{Verpoort2008}.

\subsection{Setting for first variations}\label{variation:setting_for_1st}

Let us consider a function $F$ defined on the set of immersions $\on{Imm}(M,\R^n)$, 
an immersion $f_0$ and a tangent vector $h$ to $f_0$. 
Furthermore we choose a curve of immersions 
$$f: \R \to \on{Imm}(M,\R^n) \text{ with } f(0)=f_0 \text{ and } \p_t|_0 f = h=h^{\bot}.\nu^{f_0}+Tf_0.h^{\top}.$$
In this setting, the first variation of $F$ is 
$$D_{(f_0,h)} F = T_{f_0}F.h= \p_t|_0 F(f(t)).$$

\subsection{Tangential variation of equivariant tensor fields}
\label{variation:tangential}

Let the smooth mapping $F:\on{Imm}(M,\R^n) \to \Gamma(T^r_s M)$ take values in some space of tensor fields over $M$, 
or more generally in any natural bundle over $M$, see \cite{MichorF}.
If $F$ is equivariant 
with respect to pullbacks by 
diffeomorphisms of $M$, i.e. 
$$F(f)=(\phi^* F)(f)=\phi^* (F((\phi\i)^*f)) $$ 
for all $\phi \in \on{Diff}(M)$ and $f \in \on{Imm}(M,\R^n)$,
then the tangential variation of $F$ is
\begin{align*}
D_{(f_0,Tf_0.h^\top)} F&=
\p_t|_0 F(f_0 \o Fl^{h^\top}_t)=
\p_t|_0 F\big((Fl^{h^\top}_t)^* f_0\big)\\&=
\p_t|_0 (Fl_t^{h^\top})^* (F(f_0)) = \L_{h^\top}(F(f_0)).
\end{align*}

\subsection{Variation of the volumeform \cite[section~3.5]{Michor118}}\label{variation:volume_form}

\emph{The differential of the volume form
\begin{equation*}
\left\{ \begin{array}{ccl}
\on{Imm} &\to &\Omega^{n-1}(M),\\
f &\mapsto &\on{vol}(g)=\on{vol}(f^*\bar g)
\end{array}\right.\end{equation*}
is given by}
\begin{equation*}
D_{(f_0,h)} \on{vol}(f^*\bar g) = \p_t|_0 \on{vol}(f^*\bar g) =
\big(\on{div}^{g_0}(h^\top) - \on{Tr}(L^{f_0}) .h^\bot\big) \on{vol}(g_0).
\end{equation*}

\subsection{Variation of the Weingarten map \cite[section~3.8]{Michor118}}\label{variation:weingarten}

\emph{The differential of the Weingarten map
\begin{equation*}
\left\{ \begin{array}{ccl}
\on{Imm} &\to &\Gamma(\on{End}(TM)),\\
f &\mapsto &L^f
\end{array}\right.\end{equation*}
is given by}
\begin{equation*}\begin{aligned}
D_{(f_0,h)} L^f = \p_t|_0 L&=g_0\i \nabla^2(h^\bot) + h^\bot (L^{f_0})^2 + \L_{h^\top} (L^{f_0}).
\end{aligned}\end{equation*}

\subsection{Variation of the Gau\ss curvature}\label{variation:kappa}
\emph{The differential of the Gau\ss curvature
\begin{equation*}
\left\{ \begin{array}{ccl}
\on{Imm} &\to &C^\infty(M),\\
f &\mapsto &\on{det}(L^f)
\end{array}\right.\end{equation*}
is given by}
\begin{align*}
D_{(f_0,h)} \on{det}(L^f) &= \on{Tr}(L^{f_0}). \on{det}(L^{f_0}). h^\bot
                + g^0_2\Big(g.\on{C}(L^{f_0}),\nabla^2(h^\bot) \Big)+d \on{det}(L^{f_0})(h^{\top}) 
\end{align*}
where $C(L)$ is the classical adjoint of $L$ uniquely determined by
$$\on{C}(L).L=L.\on{C}(L)=\on{det}(L).I$$
\begin{proof}
For the normal part we have 
\begin{align*}
 \p_t|_0 \on{det}(L)&= \on{Tr}\big(\on{C}(L^{f_0}) . \p_t|_0 L\big) \\
&=\on{Tr}\Big(\on{C}(L^{f_0}). \big((f_0^* \bar g)\i \nabla^2(h^\bot)
        + h^\bot (L^{f_0})^2\big)\Big)\\
&=\on{Tr}\Big(\on{C}(L^{f_0}). (f_0^* \bar g)\i .\nabla^2(h^\bot)
        + h^\bot \on{det}(L^{f_0}).L^{f_0} \Big)=\\
&=\on{Tr}(L^{f_0}). \on{det}(L^{f_0}). h^\bot + g^0_2\Big(g.\on{C}(L^{f_0}),\nabla^2(h^\bot) \Big)
\end{align*}
For the tangential part, we use the observation of section~\ref{variation:tangential}.
\end{proof}

\section{The geodesic equation on $\on{Imm}(M,\mathbb R^n)$}\label{geodesic_equation_imm}

\subsection{The geodesic equation on $\on{Imm}(M,\mathbb R^n)$}\label{geodesic_equation_imm:geodesic_equation}
We derive the geodesic equation only for $G^{\Phi}$-metrics with $\Phi=\Phi(\on{det}(L))$. 
This equation can be combined with the geodesic equation for 
metrics weighted by mean-curvature and volume \cite[section~5.1]{Michor118}.
The resulting geodesic equation for $\Phi=\Phi(\Vol,\Tr(L),\on{det}(L))$ on $\Imm$
is very long, but on $B_i$ it is much shorter and we print it in
section~\ref{geodesic_equation_Bi:geodesic_equation}.

We use the method of \cite[section~4]{Michor118} to calculate the 
geodesic equation as
\begin{align*}
&f_{tt}=\frac12 H_f(f_t,f_t)-K_f(f_t,f_t),\\
\intertext{where $H,K$ are the metric gradients given by}
&D_{f,m}G_f^\Phi(h,k)=G_f^\Phi(K_f(m,h),k)=G_f^\Phi(m,H_f(k,h)).
\end{align*}
So we need to compute the metric gradients.
The calculation also shows the existence of the gradients. 
Let $m \in T_f\on{Imm}(M,\mathbb R^n)$ with $m=m^\bot.\nu^f+Tf.m^\top.$ 
\begin{align*}
D_{(f,m)} G^\Phi_f(h,k) &=D_{(f,m)} \int_M \Phi(\det(L))\,\g(h,k)\on{vol}(g)\\              
& = \int_M \Phi'(\det(L)) (D_{(f,m)} \on{det}(L))\, \g(h,k) \on{vol}(g) \\  
& \qquad\qquad+\int_M  \Phi(\det(L))\,\g(h,k)(D_{(f,m)} \on{vol}(g)).
\end{align*}
To read off the $K$-gradient of the metric, we write this expression as
\begin{align*}
\int_M \Phi .\g\left(\big(
\frac{\Phi'}{\Phi} (D_{(f,m)} \on{det}(L))+ 
\frac{D_{(f,m)} \on{vol}(g)}{\on{vol}(g)} \big) h,k \right) \on{vol}(g)
\end{align*}
Therefore, using the formulas from section~\ref{variation} we can calculate the $K$ gradient:
\begin{align*}
K_f(m,h)&=\big[\frac{\Phi'}{\Phi} (D_{(f,m)} \on{det}(L))
+\frac{D_{(f,m)} \on{vol}(g)}{\on{vol}(g)} \big] h\\
&= \bigg[ \frac{\Phi'}{\Phi}\Big(\on{Tr}(L). \on{det}(L). m^\bot
                + g^0_2\big(g.\on{C}(L),\nabla^2(m^\bot) \big)+d\on{det}(L)(m^{\top}) \Big) \\
&\qquad+ \on{div}^g(m^\top)- \on{Tr}(L). m^\bot \bigg] h.
\end{align*}
To calculate the $H$-gradient, we treat the two summands of $D_{(f,m)} G^\Phi_f(h,k)$ separately. 
The first summand is
\begin{align*}
&\int_M \Phi' (D_{(f,m)} \on{det}(L)) \g(h,k)\on{vol}(g)=\\&\qquad
=\int_M \Phi' \on{Tr}(L). \on{det}(L). m^\bot .\g(h,k)\on{vol}(g) 
+\int_M \Phi'd\on{det}(L)(m^{\top}) \g(h,k)\on{vol}(g)\\&\qquad\qquad
+ \int_M g^0_2\big(\Phi'.\g(h,k).g.\on{C}(L),\nabla^2(m^\bot) \big)\on{vol}(g)\\&\qquad=
\int_M \Phi' \on{Tr}(L). \on{det}(L). m^\bot .\g(h,k)\on{vol}(g) \\&\qquad\qquad
+\int_M \Phi' g\big(\on{grad}^g(\on{det}(L)),m^\top\big)\g( h,k)\on{vol}(g)\\&\qquad\qquad
+ \int_M \nabla^*\nabla^*\Big(\Phi'.g.\on{C}(L)\g(h,k)\Big).m^\bot \on{vol}(g)\\&\qquad=
G^{\Phi}_f\bigg(m,\frac{\Phi'}{\Phi} \on{Tr}(L). \on{det}(L).\g(h,k).\nu+\frac{\Phi'}{\Phi}\g(h,k). Tf.\on{grad}^g(\on{det}(L))
\\&\qquad\qquad\qquad\qquad
+\frac{1}{\Phi} \nabla^*\nabla^*\big(\Phi'.g.\on{C}(L)\g(h,k)\big).\nu\bigg).
\end{align*}
In the calculation of the second term of the $H_f(m,h)$ gradient, we will make use of the following formula 
from \cite[section 5.1]{Michor118},
which is valid for $\phi\in C^\infty(M)$ and $X \in \X(M)$: 
\begin{align*}
\int_M d\phi(X)  \on{vol}(g) =-\int_M \phi . \on{div}(X)  \on{vol}(g).  
\end{align*}
Therefore we can calculate the second summand, which is given by
\begin{align*}
&\int_M \Phi .\g(h,k)(D_{(f,m)} \on{vol}(g)) 
=\int_M \Phi. \g(h,k)\big(\on{div}^g(m^\top)- \on{Tr}(L).m^\bot\big)  \on{vol}(g)\\ 
&\qquad = -\int_M d(\Phi\g(h,k))(m^\top)\on{vol}(g) +G^\Phi_f(m, -\g(h,k)\on{Tr}(L) .\nu)\\
&\qquad=  -\int_M \g\big( Tf.\on{grad}^g(\Phi\g(h,k)),m\big)\on{vol}(g) +G^\Phi_f(m, -\g(h,k)\on{Tr}(L) .\nu)\\
&\qquad=  G^\Phi_f\Big(m, -\frac{1}{\Phi} Tf.\on{grad}^g(\Phi\g(h,k))  -\g(h,k)\on{Tr}(L) .\nu\Big)\\
\end{align*}
Summing up all the terms the $H$-gradient is given by
\begin{align*}
&H_f(h,k)=\\&=\Big[\frac{\Phi'}{\Phi}.\on{Tr}(L). \on{det}(L).\g(h,k)+\frac{1}{\Phi} \nabla^*\nabla^*\big(\Phi'.g.\on{C}(L)\g(h,k)\big)
-\g(h,k)\on{Tr}(L)\Big] \nu^f \\
&\qquad 
+\frac{1}{\Phi}Tf.\Big[\Phi'\g(h,k)\on{grad}^g(\on{det}(L))
-\on{grad}^g(\Phi\g(h,k))\Big]
\end{align*}
The geodesic equation for a Gau{\ss}-curvature weighted metric $G^\Phi$ with 
$\Phi=\Phi(\on{det}(L))$ on $\on{Imm}(M,\mathbb R^n)$ is then given by
$$\boxed{\begin{aligned}
f_t &= h=h^\bot.\nu^f+Tf.h^\top,   \\
h_t&=\frac{1}{2}\Big[\frac{\Phi'}{\Phi}.\on{Tr}(L). \on{det}(L).\norm{h}^2
+\frac{1}{\Phi} \nabla^*\nabla^*\big(\Phi'.g.\on{C}(L)\norm{h}^2\big)
-\norm{h}^2 \on{Tr}(L)\Big] \nu^f\\&\qquad\qquad
+\frac{1}{2\Phi}Tf.\Big[\Phi'\norm{h}^2\on{grad}^g(\on{det}(L))- \on{grad}^g(\Phi.\norm{h}^2)\Big]\\&\qquad
-\big[\frac{\Phi'}{\Phi}\Big(\on{Tr}(L). \on{det}(L). h^\bot
+ g^0_2\big(g.\on{C}(L),\nabla^2(h^\bot) \big)+d\on{det}(L)(h^{\top}) \Big)\\&\qquad\qquad
+ \on{div}^g(h^\top)- \on{Tr}(L).h^\bot \big] h.
\end{aligned}}$$

\subsection{Momentum mappings}\label{geodesic_equation_imm:momentum_mappings}

The metric $G^\Phi$ is invariant under the action of the reparametrization group
$\on{Diff}(M)$ and under the Euclidean motion group $\mathbb R^n \rtimes
\on{SO}(n)$. According to \cite[section 4]{Michor118} the momentum mappings for these group
actions are constant along any geodesic in $\on{Imm}(M,\mathbb R^n)$:
$$\boxed{\begin{aligned}
\Phi g(f_t^\top) \on{vol}(g)\in\Ga(T^*M\otimes_M\on{vol}(M)) && \text{reparam. 
momentum}\\
\int_M \Phi f_t \on{vol}(g)\in(\mathbb R^n)^* && \text{linear momentum}\\
\int_M \Phi 
(f\wedge f_t)\on{vol}(g)\in\textstyle{\bigwedge^2}\mathbb R^n\cong \mathfrak{so}(n)^* && 
\text{angular momentum} 
\end{aligned}}$$

\section{The geodesic equation on $B_i(M,\mathbb R^n)$}\label{geodesic_equation_Bi}

\subsection{The horizontal bundle and the metric on the quotient space}\label{geodesic_equation_Bi:horizontal}

Since $\on{vol}(f^*\bar g)$, $\Tr(L)$ and $\det(L^f)$  react
equivariantly to the action of the group $\on{Diff}(M)$, every $G^\Phi$-metric is
$\on{Diff}(M)$-invariant. Thus it induces a Riemannian metric on $B_i$ such
that the projection $\pi:\on{Imm}\to B_i$ is a Riemannian submersion.
For every almost local metric $G^{\Phi}$ the horizontal bundle at the point $f$ equals the set of 
sections of the normal bundle along $f$.
Therefore the  metric on the horizontal bundle is given by 
$$ G^\Phi_f(a. \nu, b. \nu) = \int_{M} \Phi. a.b \on{vol}(f^*\bar{g}).$$ 
See \cite[section~6.1]{Michor118} for more details.

\subsection{The geodesic equation on $B_i(M,\mathbb R^n)$}\label{geodesic_equation_Bi:geodesic_equation}
The calculation of the geodesic equation can be done on the horizontal bundle instead of on $B_i$. 
This is possible because for every almost local metric a path in $B_i$ 
corresponds to exactly one horizontal path in $\on{Imm}$ (see \cite{Michor118}).
Therefore geodesics in $B_i$ correspond to horizontal geodesics in
$\on{Imm}$. A horizontal geodesic $f$ in $\on{Imm}$
has $f_t=a. \nu^f$ with $a \in C^\infty(\mathbb R \x M)$. The geodesic equation is
given by 
$$f_{tt} = \underbrace{a_t . \nu}_{\text{normal}} + \underbrace{a .
\nu_t}_{\text{tang.}} =\frac{1}{2} H(a.\nu,a.\nu)-K(a.\nu,a.\nu),$$
see \cite[section 4]{Michor118}. 
This equation splits into a normal and a tangential part. The normal part is given by
$$a_t = \g( \frac12 H(a.\nu,a.\nu)-K(a.\nu,a.\nu) , \nu ).$$
From the conservation of the reparametrization momentum, see \cite{Michor118}
it follows that the tangential part of the geodesic equation is 
satisfied automatically.

Therefore the geodesic equation on $B_i$ for $\Phi=\Phi(\on{det}(L))$ is given by 
\begin{equation*}\boxed{\begin{aligned}
f_t &=a.\nu^f,   \\
a_t&=\frac{1}{\Phi}\Big[\frac{\Phi}{2} \on{Tr}(L).a^2 
 +\frac{1}{2} \nabla^*\nabla^*\big(\Phi'.g.\on{C}(L)a^2\big)-\Phi' g^0_2\Big(g.\on{C}(L),\nabla^2(a) \Big).a\\&\qquad\qquad
-\frac{\Phi'}{2}.\on{Tr}(L). \on{det}(L).a^2\Big] .
\end{aligned}}\end{equation*}

Again the geodesic equation for a more general almost local metric $G^\Phi$
where $\Phi=\Phi(\Vol,\Tr(L),\on{det}(L))$ can be obtained by combining the above equation 
with the results in \cite[section~5.1]{Michor118}. It reads as

\begin{equation*}\boxed{\begin{aligned}
f_t &=a.\nu^f,   \\
a_t&=\frac{1}{\Phi}\Big[\frac{\Phi}{2} \on{Tr}(L).a^2 -\frac{1}{2} \on{Tr}(L) \int_M (\p_1\Phi)  a^2 \on{vol}(g)
+ (\p_1\Phi) a.\int_M \on{Tr}(L).a \on{vol}(g) \\&\qquad\qquad
-\frac{1}{2} \Delta\big((\p_2\Phi) a^2 \big) +(\p_2\Phi).a\Delta a-\frac{\p_2\Phi}{2} a^2 \on{Tr}(L^2)\\&\qquad\qquad
 +\frac{1}{2} \nabla^*\nabla^*\big((\p_3 \Phi).g.\on{C}(L)a^2\big)-(\p_3\Phi) g^0_2\Big(g.\on{C}(L),\nabla^2(a) \Big).a\\&\qquad\qquad
-\frac{\p_3\Phi}{2}.\on{Tr}(L). \on{det}(L).a^2\Big] 
\end{aligned}}\end{equation*}
For the case of curves immersed in $\R^2$, this formula specializes to the formula given in 
\cite[section~3.4]{Michor107}. (When verifying this, remember that $\Delta=-D_s^2$ and  $\on{Tr}(L)=\on{det}(L)=\ka$ in 
the notation of \cite{Michor107}.)

\section{Numerical results}\label{numerics}

We want to solve the boundary value problem for geodesics in shape space of surfaces in $\R^3$
with respect to curvature weighted metrics, 
more specifically with respect to $G^\Phi$-metrics with $\Phi=1+A.\Tr(L)^{2k}+B.\on{det}(L)^{2l}$. 
In \cite{Michor118} we did this for metrics depending on the mean curvature only.  

We will minimize horizontal path energy
\begin{equation*}
E^{\on{hor}}(f) = \int_0^1 \int_M \Phi(\Tr(L),\on{det}(L)) (f_t^\bot)^2 \vol(g)
\end{equation*}
over the set of paths $f$ of immersions with fixed endpoints. 
By definition, the horizontal path
energy does not depend on reparametrizations of the surface. Therefore
we can add a penalty term to the above energy that controls regularity of the parametrization.

We reduce this infinite dimensional problem to a finite-dimensional one 
by approximating immersed surfaces by triangular meshes. 
We solved the resulting optimization problem using the nonlinear solver IPOPT
(Interior   Point  OPTimizer  \cite{Waechter2002}).  IPOPT  was
invoked by AMPL (A Modeling Language for Mathematical Programming 
\cite{Fourer2002} ). The data file containing the definition of the combinatorics of the triangle mesh 
was automatically generated by the computer algebra system Mathematica. 

\subsection{Discrete path energy}\label{numerics:discrete}

As in \cite[section 10.1]{Michor118} we define the discrete mean curvature at vertex $p$ as
$$\Tr(L) = \frac{\norm{\text{vector mean curvature}}}{\norm{\text{vector area}}} =
\frac{\norm{\nabla (\text{surface area})}}{\norm{\nabla (\text{enclosed volume})}}.$$
Here $\nabla$ stands for a discrete gradient, and
$$(\text{vector mean curvature})_p=\big(\nabla (\text{surface area})\big)_p = 
\sum_{(p,p_i) \in E} (\cot \alpha_i + \cot \beta_i) (p-p_i)$$
is the vector mean curvature defined by the cotangent formula. In this formula,
$\alpha_i$ and $\beta_i$ are the angles opposite the edge $(p,p_i)$ in the two adjacent triangles. Furthermore, 
$$
(\text{vector area})_p=\big(\nabla (\text{enclosed volume})\big)_p 
= \sum_{f \in \on{Star}(p)} \nu(f).(\text{surface area of $f$})
$$
is the vector area at vertex $p$. 
The discrete  Gauss curvature can be defined as
$$
\on{det}(L)(p)=\frac{\ol\Theta(p)}{\text{Area of star}(p)}.
$$
Here  $\ol\Theta(p)$ stands for the angular deflection at $p$, defined by
$$\ol\Theta(p)=2\pi-\sum_{i=1}^{\#(\on{star}(p))}\theta_i,$$ 
where $\theta_i$ denotes  the internal angle of the $i$-th corner of vertex $p$ and $\#(\on{star}(p)$ the number of faces 
adjacent to vertex $p$. 
Using these definitions, the horizontal path energy and the penalty term
were discretized as in \cite[section 10.1]{Michor118}.

\subsection{Geodesics of concentric spheres}\label{numerics:concentric_conformal}
In this chapter we will study geodesics between concentric spheres for Gau{\ss} curvature weighted metrics, 
i.e. metrics with $\Phi=1+B\det(L)^{2l}$. For mean curvature weighted metrics this has been done in \cite{Michor118}.

The  set of spheres with common center $x \in \R^n$
is a totally geodesic subspace of $B_i$ with the $G^{\Phi}$-metric (see \cite[section 10.3]{Michor118}). 
Within a set of concentric spheres, any sphere is uniquely described by its radius $r$. 
Thus the geodesic equation reduces to an ordinary differential equation 
for the radius. It can be read off the geodesic equation in section~\ref{geodesic_equation_Bi}, when 
it is taken into account that all functions are constant on each sphere, and that
$$r_t=a, \quad L=-\tfrac{1}{r}\on{Id}_{TM}, \quad \Tr(L)=-\tfrac{n-1}{r}, \quad \on{det}(L)=\tfrac{1}{(-r)^{n-1}}.$$

Then the geodesic equation for $\Phi=1+B\det(L)^{2l}$ on a set of concentric spheres in $B_i$ reads as
\begin{equation*}\boxed{\begin{aligned}
r_{tt}&=-r_t^2 \frac{n-1}{2} \left( \frac{1}{r} -\frac{2l.B}{r^{(n-1)2l+1}+B.r}\right). 
\end{aligned}}\end{equation*}
Note that $\on{det}(L)=\big(\Tr(L)/(n-1)\big)^{n-1}$. 
Therefore this equation is equal to the equation for metrics 
weighted by mean curvature with suitably adapted coefficients (see \cite[section 10.3]{Michor118}).
 
{\it This equation is in accordance with the numerical results obtained by minimizing the discrete path energy. 
As will be seen, the numerics show that the shortest path connecting two concentric spheres
in fact consists of spheres with the same center, 
and that the above differential equation is satisfied. 
Furthermore, in our experiments the optimal paths obtained were independent of the initial path 
used as a starting value for the optimization. }

\begin{figure}[ht]
 \begin{psfrags}
 \def\PFGstripminus-#1{#1}%
 \def\PFGshift(#1,#2)#3{\raisebox{#2}[\height][\depth]{\hbox{%
   \ifdim#1<0pt\kern#1 #3\kern\PFGstripminus#1\else\kern#1 #3\kern-#1\fi}}}%
 \providecommand{\PFGstyle}{}%
 %
 \psfrag{B01}[cl][cl]{\PFGstyle $B=0.1$}%
 \psfrag{B100}[cl][cl]{\PFGstyle $B=100$}%
 \psfrag{B10}[cl][cl]{\PFGstyle $B=10$}%
 \psfrag{B1}[cl][cl]{\PFGstyle $B=1$}%
 \psfrag{r}[bc][bc]{\PFGstyle $r$}%
 \psfrag{t}[cl][cl]{\PFGstyle $t$}%
 \psfrag{x0}[tc][tc]{\PFGstyle $0$}%
 \psfrag{x11}[tc][tc]{\PFGstyle $1$}%
 \psfrag{x124}[tc][tc]{\PFGstyle $1200$}%
 \psfrag{x144}[tc][tc]{\PFGstyle $1400$}%
 \psfrag{x14}[tc][tc]{\PFGstyle $1000$}%
 \psfrag{x21}[tc][tc]{\PFGstyle $2$}%
 \psfrag{x23}[tc][tc]{\PFGstyle $200$}%
 \psfrag{x2}[tc][tc]{\PFGstyle $0.2$}%
 \psfrag{x41}[tc][tc]{\PFGstyle $4$}%
 \psfrag{x43}[tc][tc]{\PFGstyle $400$}%
 \psfrag{x4}[tc][tc]{\PFGstyle $0.4$}%
 \psfrag{x61}[tc][tc]{\PFGstyle $6$}%
 \psfrag{x63}[tc][tc]{\PFGstyle $600$}%
 \psfrag{x6}[tc][tc]{\PFGstyle $0.6$}%
 \psfrag{x81}[tc][tc]{\PFGstyle $8$}%
 \psfrag{x83}[tc][tc]{\PFGstyle $800$}%
 \psfrag{x8}[tc][tc]{\PFGstyle $0.8$}%
 \psfrag{y0}[cr][cr]{\PFGstyle $0$}%
 \psfrag{y11}[cr][cr]{\PFGstyle $1$}%
 \psfrag{y121}[cr][cr]{\PFGstyle $1.2$}%
 \psfrag{y141}[cr][cr]{\PFGstyle $1.4$}%
 \psfrag{y151}[cr][cr]{\PFGstyle $1.5$}%
 \psfrag{y161}[cr][cr]{\PFGstyle $1.6$}%
 \psfrag{y181}[cr][cr]{\PFGstyle $1.8$}%
 \psfrag{y1}[cr][cr]{\PFGstyle $0.1$}%
 \psfrag{y21}[cr][cr]{\PFGstyle $2$}%
 \psfrag{y2}[cr][cr]{\PFGstyle $0.2$}%
 \psfrag{y3}[cr][cr]{\PFGstyle $0.3$}%
 \psfrag{y4}[cr][cr]{\PFGstyle $0.4$}%
 \psfrag{y5}[cr][cr]{\PFGstyle $0.5$}%
 \psfrag{y6}[cr][cr]{\PFGstyle $0.6$}%
 \psfrag{ym11}[cr][cr]{\PFGstyle $-1$}%
 \psfrag{ym1253}[cr][cr]{\PFGstyle $-125$}%
 \psfrag{ym13}[cr][cr]{\PFGstyle $-100$}%
 \psfrag{ym151}[cr][cr]{\PFGstyle $-1.5$}%
 \psfrag{ym153}[cr][cr]{\PFGstyle $-150$}%
 \psfrag{ym21}[cr][cr]{\PFGstyle $-2$}%
 \psfrag{ym252}[cr][cr]{\PFGstyle $-25$}%
 \psfrag{ym52}[cr][cr]{\PFGstyle $-50$}%
 \psfrag{ym5}[cr][cr]{\PFGstyle $-0.5$}%
 \psfrag{ym752}[cr][cr]{\PFGstyle $-75$}%
 
 \includegraphics[width=\textwidth]{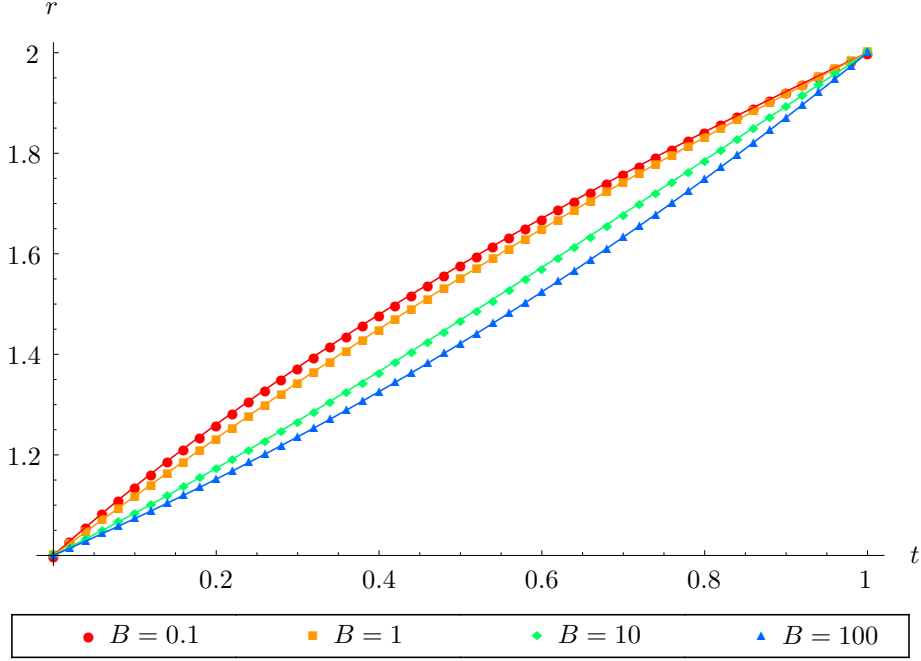}
 \end{psfrags}
 \caption{Geodesics between concentric spheres for  Gauss curvature weighted metrics with $\Phi=1+B\on{det}(L)^2$.}
 \label{fig:numerics:concentric:gb_varyingB}
 \end{figure}

A comparison of the numerical results with the exact analytic solutions can be seen in 
figure~\ref{fig:numerics:concentric:gb_varyingB}. 
The solid lines are the exact solutions. For the numerical solutions, 50 timesteps 
and a triangulation with 320 triangles were used.  

{\it The space of concentric spheres is geodesically complete with respect to the $G^{\Phi}$ metrics, 
$\Phi=1+B \on{det}(L)^{2l}$, if $l\geq 1$. }
To see this, we calculate the length of a path $f$ shrinking the unit sphere to zero ($n=3$): 
\begin{align*}
L^{G^{B,l}}_{B_i}(f) &= \int_0^1 \sqrt{G^{B,l}_f(f_t^{\bot},f_t^{\bot})} dt=
\int_0^1\sqrt{\int_M (1+B.\det(L)^{2l}) r^2_t \vol(g)}  dt\\&=
\int_0^1 |r_t| \sqrt{(1+B.\det(L)^{2l}) 4 r^2 \pi}   dt=
2\sqrt{\pi} \int_0^1 r \sqrt{(1+B.\frac{1}{r^{4l}}) }   dr.
\end{align*}
The last integral diverges if and only if $l\geq 1$.

\subsection{Translation of a sphere}

\begin{figure}[ht]
\includegraphics[width=\textwidth-10pt]{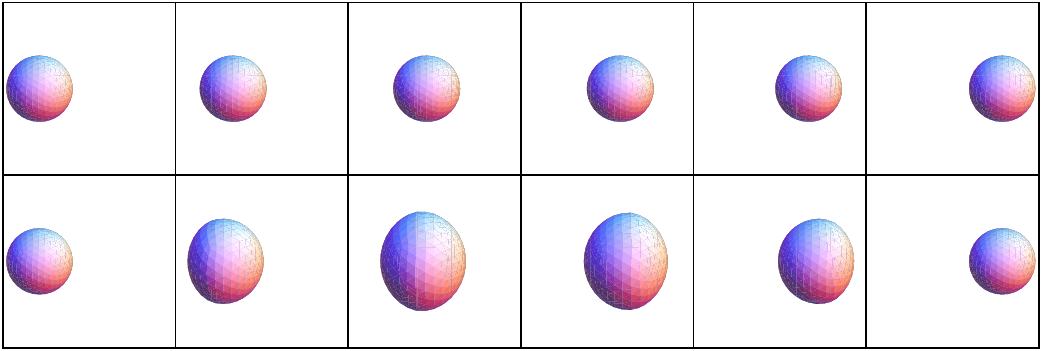}
\caption{Geodesics of translation for $\Phi=1+B\det(L)^{2}$, with $B=1$ (first row) and $B=20$ (second row) progresses from left to right.}
\label{fig:numerics:move:movegauss}
\end{figure}

For metrics with $\Phi=1+B.\det(L)^{2l}$, the following behaviours can be observed:
{\it
\begin{itemize}
\item For spheres of a certain optimal radius, pure translation is a geodesic. 
\item When the initial and final shape is a sphere with a non-optimal radius, 
the geodesic scales the sphere towards the optimal radius. 
\item When too much scaling has to be done in too little time
(always relative to the parameter $B$), 
then the geodesic passes through an ellipsoid-like shape, 
where the principal axis in the direction of the translation is shorter, 
see figure~\ref{fig:numerics:move:movegauss}.
\end{itemize}}

To determine under what conditions pure translation of a sphere is a geodesic,
let $f_t=f_0+b(t) \cdot e_1$, where $f_0$ is a sphere of radius 
$r$ and where $b(t)$ is constant on $M$. 
Plugging this into the geodesic equation~\ref{geodesic_equation_imm:geodesic_equation} 
yields an ODE for $b(t)$ 
and a part which has to vanish identically. The latter is given by:
\begin{equation*}\label{radius}
-\Phi' \frac{2}{r} \frac{1}{r^2} +\Phi  \frac{2}{r}=0 
\end{equation*}
For $\Phi=1+B.\det(L)^{2l},l\geq 1$ this yields the solution $r=\sqrt[4l]{B ( 2 l-1)}$.

\subsection{Deformation of a sphere}\label{numerics:deform}
\begin{figure}[ht]
\includegraphics[width=\textwidth-10pt]{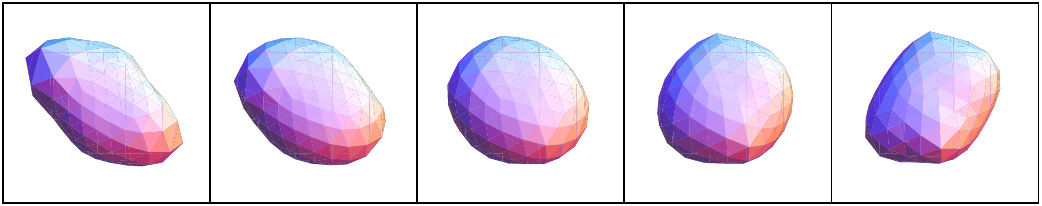}
\caption{Geodesic between two deformed shapes for $\Phi=1+\Tr(L)^2+\det(L)^{2}$. Time progresses from left to right.}
\label{fig:numerics:deform:large_deformation_gaussandmean}
\end{figure}

{\it For $\Phi=1+B \on{det}(L)$ small deformations look much the same as 
for conformal or mean curvature weighted metrics 
(compare \cite[figure~14]{Michor118}). Larger deformation are somewhat irregular,
but adding a mean curvature term yields  smooth geodesics again
(see figure \ref{fig:numerics:deform:large_deformation_gaussandmean}).
}

\bibliographystyle{plain}

\begin{thebibliography}{10}

\bibitem{Michor118}
M.~Bauer, P.~Harms, and P.~W. Michor.
\newblock Almost local metrics on shape space of hypersurfaces in n-space,
  2010.

\bibitem{Michor119}
M.~Bauer, P.~Harms, and P.~W. Michor.
\newblock Sobolev metrics on shape space of surfaces in n-space, 2010.

\bibitem{Bauer2010}
Martin Bauer.
\newblock {\em Almost local metrics on shape space of surfaces}.
\newblock PhD thesis, University of Vienna, 2010.

\bibitem{Besse2008}
Arthur~L. Besse.
\newblock {\em Einstein manifolds}.
\newblock Classics in Mathematics. Springer-Verlag, Berlin, 2008.

\bibitem{Michor40}
V.~Cervera, F.~Mascar{\'o}, and P.~W. Michor.
\newblock The action of the diffeomorphism group on the space of immersions.
\newblock {\em Differential Geom. Appl.}, 1(4):391--401, 1991.

\bibitem{Fourer2002}
R.~Fourer and B.~W. Kernighan.
\newblock {\em AMPL: A Modeling Language for Mathematical Programming}.
\newblock Duxbury Press, 2002.

\bibitem{Kobayashi1996a}
Shoshichi Kobayashi and Katsumi Nomizu.
\newblock {\em Foundations of differential geometry. {V}ol. {I}}.
\newblock Wiley Classics Library. John Wiley \& Sons Inc., New York, 1996.

\bibitem{MichorF}
I.~Kol{\'a}{\v{r}}, P.~W. Michor, and J.~Slov{\'a}k.
\newblock {\em Natural operations in differential geometry}.
\newblock Springer-Verlag, Berlin, 1993.

\bibitem{MichorG}
Andreas Kriegl and Peter~W. Michor.
\newblock {\em The convenient setting of global analysis}, volume~53 of {\em
  Mathematical Surveys and Monographs}.
\newblock American Mathematical Society, Providence, RI, 1997.

\bibitem{MichorH}
Peter~W. Michor.
\newblock {\em Topics in differential geometry}, volume~93 of {\em Graduate
  Studies in Mathematics}.
\newblock American Mathematical Society, Providence, RI, 2008.

\bibitem{Michor102}
Peter~W. Michor and David Mumford.
\newblock Vanishing geodesic distance on spaces of submanifolds and
  diffeomorphisms.
\newblock {\em Doc. Math.}, 10:217--245 (electronic), 2005.

\bibitem{Michor98}
Peter~W. Michor and David Mumford.
\newblock Riemannian geometries on spaces of plane curves.
\newblock {\em J. Eur. Math. Soc. (JEMS) 8 (2006), 1-48}, 2006.

\bibitem{Michor107}
Peter~W. Michor and David Mumford.
\newblock An overview of the {R}iemannian metrics on spaces of curves using the
  {H}amiltonian approach.
\newblock {\em Appl. Comput. Harmon. Anal.}, 23(1):74--113, 2007.

\bibitem{Verpoort2008}
Steven Verpoort.
\newblock {\em The geometry of the second fundamental form: Curvature
  properties and variational aspects}.
\newblock PhD thesis, Katholieke Universiteit Leuven, 2008.

\bibitem{Waechter2002}
A.~W\"achter.
\newblock {\em An Interior Point Algorithm for Large-Scale Nonlinear
  Optimization with Applications in Process Engineering}.
\newblock PhD thesis, Carnegie Mellon University, 2002.

\end{thebibliography}

\end{document}